\documentclass[11pt,article]{amsart}
\usepackage{amsmath, amssymb, latexsym, graphicx, mathrsfs, enumerate, amsmath, amsthm, textcomp, url, tikz-cd,bbm}
\usepackage[T1]{fontenc}
\usepackage{mathdots}
\usepackage{comment}
\usepackage[toc,page]{appendix}
\usepackage{hyperref}
\usepackage{tikz-cd}
\hypersetup{
    colorlinks=true,
    linkcolor=blue,
    filecolor=magenta,      
    urlcolor=cyan,
}

\usepackage{color}

\input xy
\xyoption{all}

\allowdisplaybreaks
\allowdisplaybreaks

\numberwithin{equation}{section}
\newtheorem{theorem}{Theorem}[section]
\newtheorem{proposition}[theorem]{Proposition}
\newtheorem{corollary}[theorem]{Corollary}
\newtheorem{lemma}[theorem]{Lemma}

\newtheorem{definition}[theorem]{Definition}

\newtheorem*{remark*}{Remark}

\newtheorem{example}[theorem]{Example}
\newtheorem{remark}[theorem]{Remark}

\newcommand{\mfo}{\mathfrak{o}}

\newcommand{\Mfo}{\mathcal{O}}

\newcommand{\cl}{\mathrm{Cl}}
\newcommand{\clb}{\overline{\mathrm{Cl}}}

\newcommand{\ord}{\mathrm{ord}}
\newcommand*{\Rom}[1]{\expandafter\@slowromancap\romannumeral #1@}
\begin{document}

\title[An upper bound for the size of the ideal class monoid]
{An upper bound for the size of the ideal class monoid}

\author[Sungmun Cho]{Sungmun Cho}
\author[Jungtaek Hong]{Jungtaek Hong}
\author[Yuchan Lee]{Yuchan Lee}
\thanks{The authors were supported by   Basic Science Research Institute Fund whose NRF grant number is RS-2021-NR060139 and the National Research Foundation of Korea(NRF) grant
funded by the Korea government(MSIT) (No. RS-2026-25508638).
}

\address{Sungmun Cho \\  Department of Mathematics, POSTECH, 77, Cheongam-ro, Nam-gu, Pohang-si, Gyeongsangbuk-do, 37673, KOREA}

\email{sungmuncho12@gmail.com}

\address{Jungtaek Hong \\  Department of Mathematics, POSTECH, 77, Cheongam-ro, Nam-gu, Pohang-si, Gyeongsangbuk-do, 37673, KOREA}

\email{jungtaekhong123@gmail.com}

\address{Yuchan Lee \\  Department of Mathematics, POSTECH, 77, Cheongam-ro, Nam-gu, Pohang-si, Gyeongsangbuk-do, 37673, KOREA}

\email{yuchanlee329@gmail.com}

\maketitle

\begin{abstract}

The ideal class monoid for an order $R$ in a finite field extension $E/F$ of a number field, denoted by $\overline{\mathrm{Cl}}(R)$, is a fundamental object to study in number theory which has useful applications in algebraic geometry and topology.  
In this paper, we describe an upper bound for $\#\overline{\mathrm{Cl}}(R)$, in terms of the class number of $E$ and (local) orbital integrals for $\mathfrak{gl}_n$. We also describe an upper bound for the class number of $E$ in terms of the \textit{Minkowski bound}.  

When $[E:F]\leq 3$ or when $R$ is a Bass order, we refine our upper bound, using a known formula for local orbital integrals in the authors' previous work. 
In particular, if $R=\mathbb{Z}[x]/(x^3-mx^2+(m-1)x-1)$ with $m\in \mathbb{Z}$ which arises in a study of \textit{Cappell-Shaneson homotopy 4-spheres} in topology, then we further refine our upper bound in terms of the discriminants of $R$ and $E$, which is $\frac{2}{3^5} \Delta_R^{\frac{1}{2}}\cdot \Delta_E^{\frac{3}{2}}$, when $\Delta_E>3075$. 
\end{abstract}

\tableofcontents

\section{Introduction}

\subsection{Background}

Let $E/F$ be a finite field extension of a number field and $R$ be an order in $E/F$.
The ideal class monoid of $R$, denoted by $\overline{\mathrm{Cl}}(R)$, is defined to be the monoid of equivalence classes of fractional $R$-ideals up to  multiplication by an element of $E^{\times}$. 
This is an important object in various areas of mathematics:
\begin{enumerate}
    \item
By \cite{Yun13}, $\#\overline{\mathrm{Cl}}(R)$ is related to the local zeta function associated with an order $R$. 
This is used to investigate the analytic properties of the global Dedekind zeta function of $R$.

    \item By \cite{Ma18}, the set of isomorphism classes of abelian varieties over $\mathbb{F}_q$ in a certain isogeny class, determined by a square-free characteristic polynomial $h(x)\in\mathbb{Q}[x]$ of Frobenius morphism (see \cite{Del69} and \cite{Cen15}),
  corresponds to the ideal class monoid of an order $\mathbb{Z}[\alpha,q/\alpha]$ in $\mathbb{Q}(\alpha)\cong\mathbb{Q}[x]/(h(x))$. Here $\mathbb{Q}(\alpha)$ is an \'etale $\mathbb{Q}$-algebra.

\item We consider a specific cubic order $R$ described as follows:
    \begin{equation}\label{eq:Caporder}
    R=\mathbb{Z}[\Theta_m]=\mathbb{Z}[x]/(x^3-mx^2+(m-1)x-1),\ m\in \mathbb{Z}.
    \end{equation}
Then determining an upper bound for $\overline{\mathrm{Cl}}(R)$ has intriguing applications in 4-dimensional manifold topology. 

The Cappell--Shaneson homotopy 4-spheres, constructed by Cappell and Shaneson in 1976 \cite{Cap761}, are parametrized by $3 \times 3$ integer matrices satisfying specific conditions. Among these, the diffeomorphism classes of smooth homotopy 4-spheres corresponding to matrices with trace $m$ are bounded from above by the cardinality of $\overline{\mathrm{Cl}}(R)$. (For further details, see \cite{AR84}, \cite{Kim23}.) 
We will discuss the properties of these cubic orders and their topological connections in greater detail in Section \ref{subsub;cs}.

\item The above equation \eqref{eq:Caporder} is also a subject for \textit{Ennola’s conjecture}. We refer to \cite{CK}.

\end{enumerate}

There has been extensive research on $\overline{\mathrm{Cl}}(R)$ in the literature, including the following:
\begin{itemize}
    \item 
\cite{DTZ} explains the minimal power $n$ such that $[I]^n$, for $[I]\in \overline{\mathrm{Cl}}(R)$, belongs to the class of an order (containing $R$) in $E$.
    \item
\cite{Ma20} explains an algorithm to compute the classes of fractional ideals in $\overline{\mathrm{Cl}}(R)$.
    \item
\cite{Ma25} presents an algorithm to compute the number of the genera in $\overline{\mathrm{Cl}}(R)$, which is equivalent to obtaining the size of $\mathrm{Cl}(R)\backslash\overline{\mathrm{Cl}}(R)$.
    \item  
\cite{Ma24} explains the ideal class monoid using the notion of \textit{Cohen-Macaulay type} of an order $R$.
When every fractional ideal of $R$ is generated by $3$ elements, Marseglia formulated $\#\overline{\mathrm{Cl}}(R)$ in terms of the class number of $R$ and the  \textit{local Cohen-Macaulay type}.
    \item
\cite{CHL} formulates an explicit closed formula for $\#(\mathrm{Cl}(R)\backslash \overline{\mathrm{Cl}}(R))$ in terms of the conductor ideal for a Bass order $R$ (see Section \ref{sec:application_bass}).
\end{itemize}

\subsection{Main Result}
In this paper, we describe an upper bound for $\#\overline{\mathrm{Cl}}(R)$ formulated as follows:
\begin{theorem}[Theorem \ref{thm:upperbound}]\label{mainthm1.2}
    We have an upper bound for $\#\overline{\mathrm{Cl}}(R)$:
    \begin{equation}\label{eq:mainthm-intro}
    \#\overline{\mathrm{Cl}}(R) \leq \#\mathrm{Cl}(\mathcal{O}_E) \prod_{\substack{v\in|\mfo|\\S_{F_v}(R_v)>0}} \#(\Lambda_{E_v}\backslash X_{R_v}),
    \end{equation}
    where $|\mfo|$ is the set of non-Archimedean places of the ring of integers $\mfo$ of $F$ with $R_v$ (resp. $E_v$) the completion of $R$ (resp. $E$) at $v$, $\#(\Lambda_{E_v}\backslash X_{R_v})$ denotes the  orbital integral associated with $R_v$ (see Definition \ref{def:orbitalintegralquo}), and $S_{F_v}(R_v)$ denotes the Serre invariant (see Definition \ref{def:generalserreinv}).
\end{theorem}

In this upper bound, $\#\mathrm{Cl}(\mathcal{O}_E)$ is the class number of a number field $E$. 
We think that an upper bound for $\#\mathrm{Cl}(\mathcal{O}_E)$ using \textit{the Minkowski bound} is widely known in the community but we were unable to find its precise statement in the literature.
Thus we describe an upper bound for $\#\mathrm{Cl}(\mathcal{O}_E)$ formulated as follows, whose proof is motivated by the finiteness theorem of the class number given in \cite[Corollary 2 of Theorem 35]{Mar18}:

\begin{theorem}[Theorem \ref{thm:clbound}]
We have an upper bound for $\#\cl(\Mfo_E)$:
    \begin{equation}\label{eq:clbound-intro}
          \#\cl(\Mfo_E)\leq \sum\limits_{\eta=1}^{\lfloor M\rfloor} \eta^{{\tilde{n}-1}}= \frac{1}{\tilde{n}}\sum\limits_{s=0}^{\tilde{n}-1}\binom{\tilde{n}}{s}B_s\lfloor M \rfloor^{\tilde{n}-s}, ~~~~~   \textit{   where}
  \end{equation}
    \begin{itemize}
        \item     $\tilde{n}=[E:\mathbb{Q}]$;
        \item $B_s$'s are the Bernoulli numbers with $B_1=\frac{1}{2}$;
        \item $M$ is the Minkowski bound $\frac{\tilde{n}!}{\tilde{n}^{\tilde{n}}} \left(\frac{4}{\pi}\right)^{r_2}|\Delta_E|^{1/2}$ (cf. \eqref{eq:minkowski}).
    \end{itemize}
    Here $2r_2$ is the number of non-real complex embeddings of $E$ and $\Delta_E$ is the discriminant of $E/\mathbb{Q}$.
   \end{theorem}

In the upper bound for $\#\overline{\mathrm{Cl}}(R)$ in Theorem \ref{mainthm1.2}, the formula for $\#(\Lambda_{E_v}\backslash X_{R_v})$ is known thanks to our previous works in \cite{JL}, \cite{CKL}, and \cite{CHL} when $[E:F]\leq 3$ or when $R$ is a Bass order. 
Thus, by applying the formula for $\#(\Lambda_{E_v}\backslash X_{R_v})$, the upper bound in Theorem \ref{mainthm1.2} becomes more explicit, as presented in Corollaries \ref{thm:n=23}-\ref{thm:bassupperbound_global1}. 

On the other hand, a main result of \cite[Theorem 6.13.(1)]{CHL} yields another upper bound for $\#\overline{\mathrm{Cl}}(R)$ using $\#\cl(R)$ (not $\#\cl(\Mfo_E)$) when $R$ is a Bass order, which is described in Corollary \ref{cor:chlthm631}. 
In Remark \ref{rmk:comparetwoex}, we compare the two upper bounds for specific examples of Bass orders and explain that their relative sharpness depends on the example.

In particular, when $R=\mathbb{Z}[x]/(x^3-mx^2+(m-1)x-1)$ as in  (\ref{eq:Caporder}), an upper bound for $\#\overline{\mathrm{Cl}}(R)$ is expressed purely in terms of the discriminants of $R$ and $E$ as follows:
\begin{corollary}[Corollary \ref{cor:main_easy}]
Let $R=\mathbb{Z}[x]/(\phi(x))$ with $\phi(x)=x^3-mx^2+(m-1)x-1$. 
If  $\Delta_E>3075$, then we have the following upper bound for $\#\overline{\mathrm{Cl}}(R)$:
\[
 \#\overline{\cl}(R) \leq \frac{2}{3^5} \Delta_\phi^{\frac{1}{2}}\cdot \Delta_E^{\frac{3}{2}}.
\]
Here $\Delta_\phi=m^4-10m^3+31m^2-30m-23$, which is the discriminant of $\phi(x)$ (and thus of $R$).
Moreover, all but finitely many $E$'s satisfy the condition that $\Delta_E>3075$ (cf. Remark \ref{rmk:final}).
\end{corollary}

We refer to Section \ref{subsub;cs}, for a background on the topological meaning of our result.

\vspace{1em}
\textbf{Organization.}
After introducing the notation in Section \ref{section:globalnotations},  we explain two local objects: orbital integrals and the relative Serre invariant in Section \ref{sec:preliminaries}.
The two main theorems are provided in Section \ref{App:AppendixB}, and Section \ref{sec:application} presents  refined upper bounds when $[E:F]\leq 3$ or when $R$ is a Bass order.
\vspace{1em}

\textbf{Acknowledgement.}
We sincerely thank Min Hoon Kim for providing a kind explanation of the background of a topological application of our results and for engaging in helpful discussions.
We are grateful to Stefano Marseglia for stimulating discussions and to Dohyeong Kim for explaining to us an upper bound for the class number of a number field.
We also sincerely thank the anonymous referee for a very careful reading of the manuscript and for numerous valuable comments and suggestions, which helped us improve the accuracy and exposition of the paper.

\section{Notations} \label{section:globalnotations}
We first introduce the notation. These are taken from \cite[Parts 1-2. Notations]{CHL}. 

\begin{itemize}
\item For a ring $A$, the set of maximal ideals is denoted by $|A|$. 
\item For a local ring $A$, the maximal ideal is denoted by $\mathfrak{m}_A$ and the residue field is denoted by $\kappa_A$. 
If $K$ is  a non-Archimedean local field, then  we sometimes use $\kappa_K$ to denote the residue field of the ring of integers in $K$, if there is no confusion. 

\item For $a\in A$ or $\psi(x) \in A[x]$ with a local ring $A$, $\overline{a}\in \kappa_A$ or $\overline{\psi}(x)\in \kappa_A[x]$ is the reduction of $a$ or $\psi(x)$ modulo $\mathfrak{m}_A$, respectively.

\item Let $F$ be a number field with $\mfo$ its ring of integers.
For $v\in |\mfo|$, 
let $F_v$ be the $v$-adic completion of $F$ with $\mfo_v$  the ring of integers of $F_v$, $\pi_v$  a uniformizer in $\mfo_v$, and $\kappa_v$ its residue field. Let $q_v=\#\kappa_v$.
For an element $x\in F_v$, $\ord_v(x)$ is the exponential valuation
with respect to $\pi_v$. 

\end{itemize}

\begin{definition}{Taken from \cite[Sections 2.1-2.3]{Ma24} or \cite[Definition A]{CHL}} \label{def:order}
Let $Z$ be a Dedekind domain with field of fractions $Q$. 
Let $K$ be an \'etale $Q$-algebra.
Note that $Z$ is always  $\mfo$ or $\mfo_v$ in this paper; in particular, all residue fields of $Z$ are finite.

\begin{enumerate}
\item{\cite[the second paragraph of Section 2.2]{Ma24}}  An order of $K$ is a subring $\Mfo$ of $K$ such that $\Mfo$ is a finitely generated $Z$-module containing $Z$ and such that $\Mfo\otimes_ZQ\cong K$. 
We say that an order $\Mfo$ is determined by an irreducible polynomial $\phi(x)\in Z[x]$ if $\Mfo \cong Z[x]/(\phi(x))$ as rings.

\item{\cite[the first paragraph of Section 2.3]{Ma24}}      A fractional $\Mfo$-ideal $I$ is a finitely generated $\Mfo$-submodule of $K$ such that $I\otimes_ZQ\cong K$.

\end{enumerate}
\end{definition}

The following notations are taken from \cite[Part 1. Notations]{CHL}.
 
\begin{itemize}
\item  The maximal order of $K$ is denoted by $\Mfo_K$. 
Here unique existence of $\Mfo_K$ is explained in \cite[the first paragraph of Section 2]{Ma20} or \cite[the third paragraph of Section 2.2]{Ma24}.
\item    The ideal quotient $(I:J)$ for two fractional $\Mfo$-ideals $I$ and $J$ is defined to be 
$    (I:J)=\{x\in K \mid xJ\subset I\}$.
Then $(I:J)$ is also a fractional $\Mfo$-ideal. 

\item A fractional $\Mfo$-ideal $I$ is called invertible if there exists a fractional $\Mfo$-ideal $J$ such that $IJ=\Mfo$. If it exists, then it is uniquely characterized by $J=\left(\Mfo:I\right)$.

\item For an order $\Mfo'$ of $K$ containing $\Mfo$, $\Mfo'$ is called an overorder of $\Mfo$.
There are only finitely many overorders of $\Mfo$ in the setting of this paper. Indeed, every overorder $\Mfo'$ is finite over $Z$, hence integral over $Z$, and therefore $\Mfo'\subset \Mfo_K$. Since $Z=\mfo$ or $\mfo_v$, all residue fields of $Z$ are finite. Hence the finite-length $Z$-module $\Mfo_K/\Mfo$ is finite as a set, and only finitely many additive subgroups of $\Mfo_K/\Mfo$ can occur as $\Mfo'/\Mfo$.

\item 
For an order $\Mfo$ of $K$, the set of fractional $\Mfo$-ideals is closed under multiplication and  the set of invertible $\Mfo$-ideals is closed under multiplication and inverse.
Thus we define the following notions:
\begin{itemize}
\item The ideal class group $\mathrm{Cl}(\Mfo)$ of $\Mfo$ is defined to be the group of equivalence classes of invertible $\Mfo$-ideals up to  multiplication by an element of $K^{\times}$. 

\item The ideal class monoid $\overline{\mathrm{Cl}}(\Mfo)$\footnote{Our convention of $\mathrm{Cl}(\Mfo)$ and $\overline{\mathrm{Cl}}(\Mfo)$ follows  \cite{Yun13}, whereas \cite{Ma24} uses $\mathrm{Pic(\Mfo)}$ and $\mathrm{ICM}(\Mfo)$ respectively.} 
of $\Mfo$ is defined to be the monoid of equivalence classes of fractional $\Mfo$-ideals up to  multiplication by an element of $K^{\times}$. 
\end{itemize}

Here, we emphasize that $(I:I)$ for $[I]\in \overline{\mathrm{Cl}}(\Mfo)$ is an order containing $\Mfo$. 
We refer to \cite[Sections 2.1-2.3]{Ma24} for detailed explanations.

\item We denote by $\Delta_{\Mfo}$ the discriminant ideal of an order $\Mfo$ over $Z$.
\begin{itemize}
\item If $\Mfo$ is determined by a polynomial $\phi(x)\in Z[x]$, then $\Delta_\Mfo$ is the same as the principal ideal $(\Delta_\phi)$ where $\Delta_\phi$ denotes the discriminant of the polynomial $\phi(x)$.
\item If $\Mfo=\Mfo_K$, then we denote $\Delta_{\Mfo}$ by $\Delta_{K/Q}$.
\item In the case that $Z$ is a PID and $v$ is a maximal ideal of $Z$, $\ord_v(\Delta_{\Mfo})$ denotes the exponential valuation of a generator of $\Delta_{\Mfo}$, as an ideal of $Z$.
\end{itemize}
\end{itemize}

From now on until the end of this paper, we fix a finite field extension $E$ of $F$ of degree $n$ and an order $R$ of $E$.
The following notation is taken from \cite[Part 2. Notations]{CHL}.

\vspace{1mm}
\[
\textit{For $v\in |\mfo|$, let  }
\left\{
\begin{array}{l}
\textit{$R_v\cong
 R\otimes_\mfo \mfo_v$ be the $v$-adic completion of $R$};\\
\textit{$E_v\cong R_v\otimes_{\mfo_v} F_v$ be the ring of total fractions of $R_v$};\\
\textit{$X_{R_v}$ be the set of fractional $R_v$-ideals so that $\overline{\cl}(R_v)=E_v^\times\backslash X_{R_v}$}.
\end{array} \right.
\]
\vspace{1mm}
\[
\textit{    For $w\in |R|$, let }
\left\{
\begin{array}{l}
\textit{$R_w$ be the $w$-adic completion of $R$};\\
\textit{$E_w\cong E\otimes_R R_w$ be the ring of total fractions of $R_w$};\\
\textit{$X_{R_w}$ be the set of fractional $R_w$-ideals  so that $\overline{\cl}(R_w)=E_w^\times\backslash X_{R_w}$}.
\end{array} \right.
\]
\vspace{1mm}

Note that $R_w$ is a local ring (possibly non-integral domain). 
By \cite[Remark B.(1)]{CHL}, we have
$R_v\cong \bigoplus\limits_{w|v,~ w\in |R|}R_w$.
Applying $(-)\otimes_{\mfo_v} F_v$ yields that
$$ E_v\cong \bigoplus\limits_{w|v,~ w\in |R|}E_w ~~~ \textit{ and }   ~~~ \Mfo_{E_v}\cong \bigoplus\limits_{w\mid v, ~ w\in |R|} \Mfo_{E_w}. $$
Thus, $E_w$ may not be a field.

\section{Local orbital integrals and the relative Serre invariant}\label{sec:preliminaries}
In this section, we will explain two local objects: orbital integrals and the relative Serre invariant. 
These two will be used to formulate an upper bound for $\#\overline{\mathrm{Cl}}(R)$.

\begin{definition}[{\cite[\S2.1 and \S4.2]{Yun13}}]\label{def:orbitalintegralquo}
Let $\mathcal{O}$ be an order of an \'etale $F_v$-algebra $K$.
Write
\[
        K=\prod_{i\in B(K)} K_i
\]
as a product of finite field extensions of $F_v$, and let
\[
        \mathcal{O}_K=\prod_{i\in B(K)}\mathcal{O}_{K_i}
\]
be the maximal order of $K$. For each $i\in B(K)$, choose a uniformizer
$\pi_i$ of $K_i$, and set
\[
        \Lambda_K=\prod_{i\in B(K)} \pi_i^{\mathbb{Z}}
        \subset \prod_{i\in B(K)}K_i^\times=K^\times .
\]
Equivalently, $\Lambda_K$ is a free abelian subgroup of $K^\times$
complementary to $\mathcal{O}_K^\times$, i.e.
\[
        K^\times=\Lambda_K\times \mathcal{O}_K^\times.
\]
Let $X_{\mathcal{O}}$ be the set of fractional $\mathcal{O}$-ideals.
The group $\Lambda_K$ acts on $X_{\mathcal{O}}$ by multiplication. We define
the orbital integral associated with $\mathcal{O}$ to be
\[
        \#(\Lambda_K\backslash X_{\mathcal{O}}).
\]
This number is independent of the choices of the uniformizers $\pi_i$.
\end{definition}\noindent

We refer the reader to \cite[Remark 2.4]{CHL} for the compatibility of this
definition with the usual orbital integrals for $\mathfrak{gl}_n$ appearing
in Arthur's trace formula. In particular, the orbital integral
$\#(\Lambda_{E_v}\backslash X_{R_v})$ is a well-defined positive integer.

\begin{definition}[{\cite[Definition 6.9]{CHL}}]\label{def:generalserreinv}
    Let $Z$ be a PID and let $Q$ be its fraction field.
    Let $K$ be an \'etale $Q$-algebra.
For an order $\Mfo$ of $K$, 
    we define the relative Serre invariant:
$$    S_Q(\Mfo)=\textit{the length of $\Mfo_K/\Mfo$ as a $Z$-module.}$$
To simplify notation, if $Q=\mathbb{Q}_p$, then we use $S_p(\Mfo)$ for $S_{\mathbb{Q}_p}(\Mfo)$.
\end{definition}

\begin{remark}\label{eq:about_serre}

\begin{enumerate}
    \item Let $F_v$ be a finite field extension of $\mathbb{Q}_p$.
If $\Mfo$ is an order of an \'etale $F_v$-algebra $K$, then $S_p(\Mfo)=[\kappa_v:\mathbb{F}_p]\cdot S_{F_v}(\Mfo)$.

\item If $\mathcal{O}$ is an order of an \'etale $F_v$-algebra $K$, then 
\begin{equation}\label{eq:serreinv}
S_{F_v}(\Mfo)=\frac{1}{2}(\ord_v(\Delta_\Mfo)-\ord_v(\Delta_{K/F_v})),
\end{equation}
which follows from the proof of \cite[Proposition 2.5]{CKL}.

\item 
Recall that $R$ is an order of a finite field extension $E$ of a number field $F$.
Then, since $\Delta_{R}$ and $\Delta_{E/F}$ are ideals in $\mfo$, both $\ord_v(\Delta_R)$ and $\ord_v(\Delta_{E/F})$ are zero for all but finitely many $v\in|\mfo|$.
Since the discriminant and the completion commute with each other by \cite[Chapter I.3, Proposition 6.(iii)]{Cas}, 
both $\ord_v(\Delta_{R_v})$ and $\ord_v(\Delta_{E_v/F_v})$ are also zero for all but finitely many $v\in|\mfo|$.
Thus for such $v$, (\ref{eq:serreinv}) gives that $S_{F_v}(R_v)=0$, which means that $R_v$ is the ring of integers in $E_v$. 
In this case, $R_v=\mathcal{O}_{E_v}$ is a product of discrete valuation rings. Hence every fractional $R_v$-ideal is of the form $\lambda R_v$ for some $\lambda\in\Lambda_{E_v}$, so the action of $\Lambda_{E_v}$ on $X_{R_v}$ is transitive and the orbital integral associated with $R_v$ is $1$.

\end{enumerate}
\end{remark}

\section{An upper bound for $\#\overline{\mathrm{Cl}}(R)$} \label{App:AppendixB}
A goal of this section is to provide an upper bound for $\#\overline{\mathrm{Cl}}(R)$ in Theorem \ref{thm:upperbound}.
We begin with the following identity, which is a reformulation of \cite[Equation (3.4)]{Yun13}.
\begin{lemma}\label{thm_in_yun}
We have the following equation
\begin{equation}\label{eq:yunyun}
\#\cl(\Mfo_E)\prod\limits_{v\in |\mfo|}\#\left(\Lambda_{E_v}\backslash X_{R_v}  \right)=\#\cl(R)\sum\limits_{\{I\}\in \cl(R)\backslash \clb(R)}\frac{\#\left(\Mfo_E^{\times}/\mathrm{Aut}(I)  \right)}{\#\mathrm{Stab}_{\cl(R)}([I])}.
\end{equation}
Here, for a fractional $R$-ideal $I$, we use the following notations:
\[\left\{
\begin{array}{l}
    \textit{$\{I\}$ is the $\cl(R)$-orbit of $[I]\in \overline{\mathrm{Cl}}(R)$};\\
     \textit{$\mathrm{Stab}_{\cl(R)}([I])$ denotes the stabilizer of $[I]$ under the $\mathrm{Cl}(R)$-action on $\overline{\mathrm{Cl}}(R)$};\\
      \textit{$\mathrm{Aut}(I)=\{x\in E^{\times}\mid xI=I\}$}.
\end{array}
\right.
\]

\end{lemma}

\begin{proof}
We start with \cite[Equation (3.4)]{Yun13} stated below:
\begin{equation}\label{eq:yun}
\#\cl(\Mfo_E)\prod\limits_{w\in|R|}\#\left(\Lambda_{E_w}\backslash X_{R_w}  \right)=\#\cl(R)\sum\limits_{\{I\}\in \cl(R)\backslash \clb(R)}\frac{\#\left(\Mfo_E^{\times}/\mathrm{Aut}(I)  \right)}{\#\mathrm{Stab}_{\cl(R)}([I])}.
\end{equation}
Thus it suffices to show that 
\begin{equation}\label{eq:relbetvnw}
\prod\limits_{w\in|R|}\#\left(\Lambda_{E_w}\backslash X_{R_w}  \right)=
\prod_{v\in|\mfo|}\prod\limits_{\substack{w\in|R|, \\ w|v}}\#\left(\Lambda_{E_w}\backslash X_{R_w}  \right)=
\prod\limits_{v\in|\mfo|}\#\left(\Lambda_{E_v}\backslash X_{R_v}  \right),
\end{equation}

where $\Lambda_{E_w}$ is defined as in
Definition~\ref{def:orbitalintegralquo}; equivalently, if
$E_w=\prod_i E_{w,i}$, then
$\Lambda_{E_w}=\prod_i \pi_{E_{w,i}}^{\mathbb{Z}}$.

The first identity in (\ref{eq:relbetvnw}) follows from the surjectivity of the map from $\mathrm{Spec}~R$ to $\mathrm{Spec}~\mfo$. 
For the second identity in (\ref{eq:relbetvnw}), the isomorphism
$R_v\cong \bigoplus\limits_{w\in |R|,\ w|v}R_w$ in
\cite[Remark B.(1)]{CHL} implies that
$X_{R_v}\cong \prod_{w|v}X_{R_w}$.
Under the corresponding decomposition
$E_v\simeq \prod_{w|v}E_w$, the choice of $\Lambda$ in
Definition~\ref{def:orbitalintegralquo} is compatible with products, namely
$\Lambda_{E_v}\cong \prod_{w|v}\Lambda_{E_w}$. Hence
\[
        \Lambda_{E_v}\backslash X_{R_v}
        \cong
        \prod_{\substack{w\in|R|, \\ w|v}}(\Lambda_{E_w}\backslash X_{R_w}).
\]
\end{proof}

To define stratifications of $\overline{\mathrm{Cl}}(R)$ and $\mathrm{Cl}(R)\backslash \overline{\mathrm{Cl}}(R)$ from \cite[Section 5]{CHL}, we introduce the following notion.

\begin{definition}[{\cite[Definition 5.1]{CHL}}]\label{def:globalclcl}
For an overorder $\mathcal{O}$ of $R$, we define the following sets:
\[\left\{
\begin{array}{l}
     \mathrm{cl}(\mathcal{O})=\{ [I] \in \overline{\mathrm{Cl}}(R) \mid (I:I)=\mathcal{O} \};  \\
     \overline{\mathrm{cl}(\mathcal{O})}=\{\textit{the set of $\mathrm{Cl}(R)$-orbits in $\mathrm{cl}(\mathcal{O})$}\}=\{ \{I\}\in\mathrm{Cl}(R)\backslash \overline{\mathrm{Cl}}(R) \mid (I:I)=\mathcal{O}  \}.
\end{array}\right.
\]    
\end{definition}
Here the set  $\overline{\mathrm{cl}(\mathcal{O})}$ is well-defined since $(JI:JI)=JJ^{-1}(I:I)=(I:I)$ for $J\in \mathrm{Cl}(R)$.
Note that $\mathrm{cl}(\mathcal{O})$ is denoted by $\mathrm{ICM}_{\Mfo}(R)$ in \cite{Ma24}.
Since $(\Mfo:\Mfo)=\Mfo$, two different overorders of $R$ are in different $\mathrm{Cl}(R)$-orbits.

\begin{proposition}[{\cite[Proposition 5.2]{CHL}}]\label{prop:stratforglobal}
We have the following stratifications:
\[
\overline{\mathrm{Cl}}(R)=\bigsqcup_{R\subset \mathcal{O} \subset \mathcal{O}_E}\mathrm{cl}(\mathcal{O})
~~~~~~~~~\textit{ and } ~~~~~~~~~~~  
\mathrm{Cl}(R)\backslash \overline{\mathrm{Cl}}(R)
=\bigsqcup_{R\subset \mathcal{O}\subset \mathcal{O}_E}\overline{\mathrm{cl}(\mathcal{O})},
\]  
where the index runs over all overorders $\mathcal{O}$ of $R$.
\end{proposition}

\begin{proposition}\label{prop:maineqofglobal}
We have the following equation
\[    
    \prod\limits_{v\in |\mfo|}\#\left(\Lambda_{E_v}\backslash X_{R_v}  \right)
    =\frac{1}{\#\mathrm{Cl}(\mathcal{O}_E)}\sum\limits_{R\subset \mathcal{O}\subset \mathcal{O}_E}\#\mathrm{cl}(\mathcal{O})\cdot \#(\mathcal{O}_E^\times/\mathcal{O}^{\times}),
\]    
where the sum is taken over all overorders $\mathcal{O}$ of $R$.
\end{proposition}

\begin{proof}
We rewrite Lemma \ref{thm_in_yun} as follows:
\begin{equation}\label{eq:intermidiate}
\prod_{v\in|\mfo|}\#\bigl(\Lambda_{E_v}\backslash X_{R_v}\bigr)
= \frac{1}{\#\mathrm{Cl}(\mathcal{O}_E)}
\sum_{\{I\}\in \cl(R)\backslash \clb(R)}
\#\bigl(\mathcal{O}_E^{\times}/\mathrm{Aut}(I)\bigr)\,
\frac{\#\mathrm{Cl}(R)}{\#\mathrm{Stab}_{\mathrm{Cl}(R)}([I])}.
\end{equation}
Note that $\mathrm{Aut}(I)=(I:I)^\times$.

By the stratification
$\mathrm{Cl}(R)\backslash \overline{\mathrm{Cl}}(R)
= \bigsqcup\limits_{R\subset \mathcal{O}\subset \mathcal{O}_E}
\overline{\mathrm{cl}(\mathcal{O})}$ given in Proposition \ref{prop:stratforglobal}, we rewrite the index set of the summation on the right-hand side of (\ref{eq:intermidiate}) using the overorders of $R$ so that we have
\begin{align*}
\prod_{v\in|\mfo|}\#\bigl(\Lambda_{E_v}\backslash X_{R_v}\bigr)
&= \frac{1}{\#\mathrm{Cl}(\mathcal{O}_E)}
\sum_{R\subset \mathcal{O}\subset \mathcal{O}_E}
\Biggl(
\#\bigl(\mathcal{O}_E^{\times}/\mathcal{O}^{\times}\bigr)
\sum_{\{I\}\in \overline{\mathrm{cl}(\mathcal{O})}}
\frac{\#\mathrm{Cl}(R)}{\#\mathrm{Stab}_{\mathrm{Cl}(R)}([I])}
\Biggr) .
\end{align*}

We observe that each summand $\frac{\#\mathrm{Cl}(R)}{\#\mathrm{Stab}_{\mathrm{Cl}(R)}([I])}$
is the cardinality of the $\mathrm{Cl}(R)$-orbit of $[I]$ in $\mathrm{cl}(\mathcal{O})$. 
Therefore, the sum of $\frac{\#\mathrm{Cl}(R)}{\#\mathrm{Stab}_{\mathrm{Cl}(R)}([I])}$ over $\{I\}$ in $\overline{\mathrm{cl}(\mathcal{O})}$ is exactly the same as $\#\mathrm{cl}(\mathcal{O})$. This completes the proof. 
\end{proof}

\begin{theorem}\label{thm:upperbound}
    We have an upper bound for $\#\overline{\mathrm{Cl}}(R)$:
    \begin{equation}\label{eq:mainthm}
    \#\overline{\mathrm{Cl}}(R) \leq \#\mathrm{Cl}(\mathcal{O}_E) \prod_{\substack{v\in|\mfo|\\S_{F_v}(R_v)>0}} \#(\Lambda_{E_v}\backslash X_{R_v})
    =\#\mathrm{Cl}(\mathcal{O}_E) \prod_{\substack{w\in|R|\\ S_{F_v}(R_w)>0}} \#(\Lambda_{E_w}\backslash X_{R_w}).
    \end{equation}
\end{theorem}
\begin{proof}
    Since $\#\overline{\mathrm{Cl}}(R)=\sum\limits_{R\subset \mathcal{O}\subset \Mfo_E}\#\mathrm{cl}(\mathcal{O})$ by Proposition \ref{prop:stratforglobal}, we have
    \[
    \#\overline{\mathrm{Cl}}(R)=\sum_{R\subset \mathcal{O}\subset\mathcal{O}_E}\#\mathrm{cl}(\mathcal{O})\leq
    \sum\limits_{R\subset \mathcal{O}\subset \mathcal{O}_E}\#\mathrm{cl}(\mathcal{O})\cdot \#(\mathcal{O}_E^\times/\mathcal{O}^{\times})
    =\#\mathrm{Cl}(\mathcal{O}_E)\prod\limits_{v\in |\mfo|}\#\left(\Lambda_{E_v}\backslash X_{R_v}  \right).
    \]    
Here the inequality follows from the fact that  $\#(\mathcal{O}_E^{\times}/\mathcal{O}^\times)\geq 1$ and the second identity follows from Proposition \ref{prop:maineqofglobal}.
Note that the product is taken over all $v \in |\mfo|$ but, by Remark \ref{eq:about_serre}.(2), only those finitely many places with $S_{F_v}(R_v) > 0$ contribute non-trivially.
    On the other hand, the isomorphism $\Lambda_{E_v}\backslash X_{R_v} \cong  \prod_{w|v}\Lambda_{E_w}\backslash X_{R_w}$ as described in the proof of Lemma \ref{thm_in_yun} yields the second equality in (\ref{eq:mainthm}).
    This completes the proof.
\end{proof}

We now compute an upper bound for $\#\cl(\Mfo_E)$ using Minkowski's bound.
By \cite[Corollary 2 of Theorem 37]{Mar18}, there exists a set of representatives for the ideal class group $\mathrm{Cl}(\mathcal{O}_E)$ consisting of integral ideals $I$ such that 
\begin{equation}
\label{eq:minkowski}
\mathrm{Norm}(I)\leq \frac{\tilde{n}!}{\tilde{n}^{\tilde{n}}} \left(\frac{4}{\pi}\right)^{r_2}|\Delta_E|^{1/2}.
\end{equation}
Here $2r_2$ is the number of non-real complex embeddings of $E$ and $\tilde{n}=[E:\mathbb{Q}]$. 
This upper bound is called the Minkowski bound.

\begin{theorem} \label{thm:clbound}
We have the following upper bound for $\#\cl(\Mfo_E)$:
    \begin{equation}\label{eq:clbound}
        \#\cl(\Mfo_E)\leq \sum\limits_{\eta=1}^{\lfloor M\rfloor} \eta^{{\tilde{n}-1}}= \frac{1}{\tilde{n}}\sum\limits_{s=0}^{\tilde{n}-1}\binom{\tilde{n}}{s}B_s\lfloor M \rfloor^{\tilde{n}-s}.
    \end{equation}
    Here, $B_s$'s are the Bernoulli numbers with $B_1=\frac{1}{2}$, and $M$ is the Minkowski bound $\frac{\tilde{n}!}{\tilde{n}^{\tilde{n}}} \left(\frac{4}{\pi}\right)^{r_2}|\Delta_E|^{1/2}$ in \eqref{eq:minkowski}.
\end{theorem}

\begin{proof}
The proof is motivated by the argument for the finiteness of the class number in \cite[Corollary 2 of Theorem 35]{Mar18}.
We fix a set of ideals $\mathcal{I}(\mathcal{O}_E)$ in $\mathcal{O}_E$ representing all classes in the ideal class group $\mathrm{Cl}(\mathcal{O}_E)$ satisfying the inequality (\ref{eq:minkowski}), so that $\#\mathcal{I}(\mathcal{O}_E)=\#\mathrm{Cl}(\mathcal{O}_E)$.
We give an upper bound for the number of ideals $I$ in $\mathcal{I}(\Mfo_E)$ with fixed norm $\eta$. 
Let $\eta=p_1^{k_1}\cdots p_r^{k_r}$ be the prime factorization of $\eta$.
Since $\Mfo_E$ is a Dedekind domain, $I$ factors uniquely as
\[ I=I_1\cdots I_r, \]
where each $I_i$ is an ideal whose prime factors lie over $(p_i)$, and whose norm is $p_i^{k_i}$.
For each $i$, the number of choices of $I_i$ is maximized when $p_i\Mfo_E$ splits into $\tilde{n}$ prime ideals having norm $p_i$.
In this case, the number of such $I_i$'s equals the number of non-negative integer solutions $(x_1,\cdots,x_{\tilde{n}})$ to $x_1+\cdots+x_{\tilde{n}}=k_i$, which is bounded above by $\tilde{n}^{k_i}$.
Therefore, the total number of ideals in $\mathcal{I}(\Mfo_E)$ with norm $\eta$ is at most $\tilde{n}^{k_1+\cdots+k_r}$.

Since $2^{k_1+\cdots+k_r}\leq \eta$, we have $\tilde{n}^{k_1+\cdots+k_r}\leq \eta^{\log_2{\tilde{n}}}\leq \eta^{\tilde{n}-1}$.
This yields that \[\#\cl(\Mfo_E) \leq \sum\limits_{\eta=1}^{\lfloor M \rfloor} \eta^{\tilde{n}-1},\]
where $M$ is the Minkowski bound $\frac{\tilde{n}!}{\tilde{n}^{\tilde{n}}} \left(\frac{4}{\pi}\right)^{r_2}|\Delta_E|^{1/2}$ in \eqref{eq:minkowski}.
This yields the desired bound.
The second equality of \eqref{eq:clbound} follows from Faulhaber's formula.
\end{proof}

\section{Refined upper bound from explicit value for orbital integrals}\label{sec:application}
In this section, we explain a refined upper bound for $\#\overline{\mathrm{Cl}}(R)$ in the following cases, by explicit value of orbital integrals:
\begin{enumerate}
    \item $R\cong \mfo[x]/(\phi(x))$ where $\phi(x)$ is irreducible with  $\deg(\phi)=2$ or $3$;
    \item $R$ is a Bass order;
    \item $R\cong \mathbb{Z}[x]/(x^3-mx^2+(m-1)x-1)$.
\end{enumerate}
The third case is a special case of the first  and the second case. 
Nonetheless we treat it independently due to its application in counting similarity classes of Cappell–Shaneson matrices (see Section \ref{subsec:Cap}).

\subsection{The case that $R\cong \mfo[x]/(\phi(x))$ where $\phi(x)$ is irreducible with  $\deg(\phi(x))=2$ or $3$}\label{sec:application23}
We set up notations as follows: 
\begin{itemize}
    
    \item For a polynomial $\phi(x)$, $B_v(\phi)$ is an index set in bijection with irreducible factors $\phi_{v,i}$ of $\phi$ over $F_v$, so that
    \[
    \left\{
    \begin{array}{l}
         R_v\cong \bigoplus_{i\in B_v(\phi)}R_{v,i}\textit{ where $R_{v,i}\cong \mfo_v[x]/(\phi_{v,i}(x))$};\\
         E_v\cong \bigoplus_{i\in B_v(\phi)}E_{v,i}\textit{ where $E_{v,i}\cong F_v[x]/(\phi_{v,i}(x))$}.
    \end{array}
    \right.
    \]

    \item We set $T_{E_v}(\kappa_v)=(\mathcal{O}_{E_v}\otimes_{\mfo_v}\kappa_v)^\times$ so that
   \[ \#T_{E_v}(\kappa_v)=\prod_{i\in B_v(\phi)}q_v^{[E_{v,i}:F_v]}(1-\frac{1}{\#\kappa_{E_{v,i}}}),
   \]
    where $\kappa_{E_{v,i}}$ is the residue field of $E_{v,i}$.
\end{itemize}

\cite[Proposition 7.3]{JL} gives an explicit value for the product of $\#(\Lambda_{E_v}\backslash X_{R_v})$ appearing in Theorem \ref{thm:upperbound}. Note that this value is based on the closed formula for $\#(\Lambda_{E_v}\backslash X_{R_v})$ with $v \in |\mfo|$ provided in \cite[Sections 5-6]{CKL}.
Thus we have the following result:

\begin{corollary}\label{thm:n=23}
\begin{enumerate}
    \item If $\deg(\phi(x))=2$, then we have
\begin{equation}\label{thmn=2}
\#\overline{\mathrm{Cl}}(R)\leq \#\mathrm{Cl}(\mathcal{O}_E) \prod_{\substack{v\in|\mfo|\\S_{F_v}(R_v)>0}}\left( 1+\frac{\#T_{E_v}(\kappa_v)}{q_v-1}\cdot \frac{q_v^{S_{F_v}(R_v)}-1}{q_v-1} \right),
\end{equation}
\item If $\deg(\phi(x))=3$, then we have
\begin{equation}\label{eq:thmn=3}
    \#\overline{\mathrm{Cl}}(R)\leq \#\mathrm{Cl}(\mathcal{O}_E)\prod_{\substack{v\in|\mfo|\\S_{F_v}(R_v)>0}}q_v^{\rho_v(\phi)}\left(  1+\frac{\#T_{E_v}(\kappa_v)}{(q_v-1)^2}\cdot\mathcal{F}_{v}(\phi)\right),
    \end{equation}
Here we use the following notations for $v\in|\mfo|$:
\begin{align*}
\rho_v(\phi) &=
\sum_{\{i,j\}\subset B_v(\phi),\,i\neq j} \ord_v(\mathrm{Res}(\phi_{v,i},\phi_{v,j})) =
\begin{cases}
0 & \textit{if } \#B_v(\phi)=1; \\
S_{F_v}(R_v)-\delta_v & \textit{if } \#B_v(\phi)=2; \\
S_{F_v}(R_v) & \textit{if } \#B_v(\phi)=3.
\end{cases} \\
\\
\mathcal{F}_v(\phi) &=
\begin{cases}
\dfrac{q_v^{\delta_v}-1}{q_v-1}-\dfrac{3(q_v^{\delta_v-d_v}-1)}{q_v^2-1} 
    & \textit{if $E_v/F_v$ is unramified}; \\[1ex]
\dfrac{q_v^{\delta_v}-1}{q_v-1}-\dfrac{3(q_v^{\delta_v-d_v}-1)}{q_v^2-1}
 + \dfrac{(1+\delta_v-3d_v)q_v^{\delta_v-d_v}-1}{q_v(q_v+1)}
    & \textit{if $E_v/F_v$ is ramified}; \\[2ex]
\dfrac{q_v^{\delta_v}-1}{q_v-1}
    & \textit{if $\phi(x)$ is reducible over $F_v$}.
\end{cases}
\end{align*}
where $\mathrm{Res}(\cdot,\cdot)$ denotes the resultant of two polynomials, $\delta_v=\max\limits_{i \in B_v(\phi)}\{S_{F_v}(R_{v,i})\}$, and $\ d_v=\lfloor \frac{\delta_v}{3} \rfloor$.
\end{enumerate}
In the above two cases, $S_{F_v}(R_v)>0$ if and only if $R_v$ is not the maximal order in $F_v[x]/(\phi(x))$. 
\end{corollary}
See Theorem \ref{thm:clbound} for an upper bound for  $\#\cl(\Mfo_E)$.

\subsection{The case that $R$ is a Bass order}\label{sec:application_bass}
A Bass order $R$ is an order of $E$ such that every ideal of $R$ is generated by $2$ elements.

When $R$ is a Bass order, a closed formula for $\#(\Lambda_{E_w}\backslash X_{R_w})$ with $w\in |R|$ is given in \cite{CHL}.
Based on this formula, we derive an upper bound for $\#\overline{\mathrm{Cl}}(R)$ in terms of the class number of $E$.
We set up notations as follows:
\begin{itemize}
    \item Recall that $\kappa_{R_w}$ denotes the residue field of a local ring $R_w$. We denote $\#\kappa_{R_w}$ by $q_{R_w}$. 
    \item For $w\in |R|$ and $v\in |\mfo|$ with $w \mid v$, we denote  by $K_w$ the unramified field extension of $F_v$ in $E_w$ corresponding to the field extension $\kappa_{R_w} / \kappa_v$ (see \cite[Section 1.2]{CHL}). 

    It follows that $q_{R_w}^{S_{K_w}(R_w)}=q_v^{S_{F_v}(R_w)}$ since $S_{K_w}(R_w)\cdot [\kappa_{R_w}:\kappa_v]=S_{F_v}(R_w)$ due to Remark \ref{eq:about_serre}.(1).
    Here, the proof of \cite[Theorem 3.1]{CKL} yields that $\mathcal{O}_{K_w}\subset R_w$.
    \item We write 
$$|R|=|R|^{irred}\sqcup |R|^{split} ~~~ \textit{ where  }  ~~~ \left\{
      \begin{array}{l}
|R|^{irred}\subset  \{w\in |R| : \textit{ $R_w$ is an integral domain}\};\\
|R|^{split}\subset  \{w\in |R| : \textit{ $R_w$ is not an integral domain}\}.
      \end{array} \right.$$ 
\end{itemize}

\begin{corollary}\label{thm:bassupperbound_global1}
    For a Bass order $R$ of a number field $E$, we have
\[
\#\overline{\mathrm{Cl}}(R)\leq \# \mathrm{Cl}(\mathcal{O}_E) \left(\prod_{\substack{w\in |R|^{split}\\ S_{K_w}(R_w)>0}}q_{R_w}^{S_{K_w}(R_w)}  \prod_{\substack{w\in |R|^{irred}\\ S_{K_w}(R_w)>0}}
\Bigg(
q_{R_w}^{S_{K_w}(R_w) }+[\kappa_{E_w}:\kappa_{R_w}]\cdot \frac{q_{R_w}^{S_{K_w}(R_w)}-1}{q_{R_w}-1}\Bigg)
\right).
\]
\end{corollary}
\begin{proof}    
    Theorem \ref{thm:upperbound} provides the inequality
    \[
    \#\overline{\mathrm{Cl}}(R)\leq \#\mathrm{Cl}(\mathcal{O}_E)\prod_{\substack{w\in |R|\\ S_{K_w}(R_w)>0}}\#(\Lambda_{E_w}\backslash X_{R_w}).
    \]
    By Proposition \cite[Proposition 6.3.(2)]{CHL}, $R_w$ is a Bass order of $E_w$.
    Then the computation of $\#(\Lambda_{E_w}\backslash X_{R_w})$ for each $w\in |R|$ proceeds as follows.
    \begin{itemize}
        \item 
    For $w\in |R|^{irred}$, by \cite[Theorem 3.7]{CHL}, we have
    \[
    \#(\Lambda_{E_w}\backslash X_{R_w})=
    q_{R_w}^{S_{K_w}(R_w)}+[\kappa_{E_w}:\kappa_{R_w}]\cdot (q_{R_w}^{S_{K_w}(R_w)-1}+q_{R_w}^{S_{K_w}(R_w)-2}+\cdots+ q_{R_w}+1).
    \]
        \item 
    For $w\in |R|^{split}$, by \cite[Theorem 6.11.(1)]{CHL}, we have 
    \[\#(\Lambda_{E_w}\backslash X_{R_w})=q_{R_w}^{S_{K_w}(R_w)}.\]
    \end{itemize}
\end{proof}
See Theorem \ref{thm:clbound} for an upper bound for  $\#\cl(\Mfo_E)$.

A main result of \cite{CHL} also yields another upper bound for $\#\overline{\mathrm{Cl}}(R)$. We will state it in the following corollary.
A major difference between the two upper bounds is:
\begin{itemize}
    \item the upper bound in Corollary \ref{thm:bassupperbound_global1} involves $\#\mathrm{Cl}(\Mfo_E)$ whose upper bound is described in Theorem \ref{thm:clbound};
    \item the upper bound in Corollary \ref{cor:chlthm631} involves $\#\mathrm{Cl}(R)$ whose upper bound is unknown.
    \end{itemize}

\begin{corollary}\label{cor:chlthm631}
For a Bass order $R$, we have 
\[
 \#\overline{\mathrm{Cl}}(R)\leq \#\mathrm{Cl}(R)\prod_{w|\mathfrak{f}(R)}\left( S_{K_w}(R_w)+1 \right).
\]
Here $\mathfrak{f}(R)$ is the conductor ideal of $R$, which is defined to be the biggest ideal of $\Mfo_E$ which is contained in $R$. 
\end{corollary}
\begin{proof}

For a Bass order $R$, by \cite[Proposition 6.3.(1) and Theorem 6.13.(1)]{CHL}, we have
\[
\#\overline{\mathrm{Cl}}(R)=\sum\limits_{R\subset \Mfo \subset \Mfo_E}\#\cl(\Mfo) ~~~ \textit{  and   }  ~~~~
\#(\mathrm{Cl}(R)\backslash \overline{\mathrm{Cl}}(R))=\prod\limits_{w|\mathfrak{f}(R)}\left( S_{K_w}(R_w)+1 \right).\]

For any overorder $\mathcal{O}$ of $R$, \cite[Remark 3.8]{Ma20} shows that $\#\mathrm{Cl}(\mathcal{O}) \leq \#\mathrm{Cl}(R)$.
In addition, \cite[Proposition 6.3.(1)]{CHL} yields that the number of overorders of $R$ equals $\#(\mathrm{Cl}(R)\backslash \overline{\mathrm{Cl}}(R))$.
Thus we have the desired upper bound. 
\end{proof}

Note that the explicit description of the conductor ideal $\mathfrak{f}(R)$ is explained in \cite[Section 6.3]{CHL}. 

\begin{example}\label{ex:bounds}
    We compare the two upper bounds in Corollaries \ref{thm:bassupperbound_global1}-\ref{cor:chlthm631} for two explicit examples of quadratic orders. Note that any quadratic order is a Bass order.
     \begin{enumerate}
        \item Let $R=\mathbb{Z}[3\sqrt{2}]\subset E=\mathbb{Q}[\sqrt{2}]$. 
        Then the maximal order is $\Mfo_E=\mathbb{Z}[\sqrt{2}]$, which is a PID, so that $\#\cl(\Mfo_E)=1$. 
        The Serre invariant of $R_w$ at each place $w\in |R|$ is $$S_{K_w}(R_w)=\left\{\begin{array}{l l}
1& \textit{if $w=(3)$};\\
0&  \textit{otherwise.}
\end{array}\right.$$
For $w=(3)$, the local ring $R_w$ is an integral domain with $[\kappa_{E_w}:\kappa_{R_w}]=2$ and $q_{R_w}=3$. 
Thus the upper bound in Corollary \ref{thm:bassupperbound_global1} is $5$.

On the other hand, a direct computation using Magma  \cite{BCP97} shows that $\#\cl(R)=1$. 
Since $\#(\Mfo_E/R)=3$, the ideal $(3)$ divides the conductor ideal $\mathfrak{f}(R)$.
Thus the upper bound in Corollary \ref{cor:chlthm631} is $2$.

Using Magma, we see that $\#\clb(R)=2$, so that the second upper bound is sharp.
\item
    Let $R=\mathbb{Z}[9\sqrt{-1}]\subset E=\mathbb{Q}[\sqrt{-1}]$.
    Then the maximal order is $\Mfo_E=\mathbb{Z}[\sqrt{-1}]$, which is a PID, so that $\#\cl(\Mfo_E)=1$. 
    The Serre invariant of $R_w$ at each place $w\in |R|$ is $$S_{K_w}(R_w)=\left\{\begin{array}{l l}
2& \textit{if $w=(3)$};\\
0&  \textit{otherwise.}
\end{array}\right.$$
For $w=(3)$, the local ring $R_w$ is an integral domain with $[\kappa_{E_w}:\kappa_{R_w}]=2$ and $q_{R_w}=3$. 
Thus the upper bound in Corollary \ref{thm:bassupperbound_global1} is $17$.

On the other hand, a direct computation using Magma  \cite{BCP97} shows that $\#\cl(R)=6$.
Since $\#(\Mfo_E/R)=9$, the ideal $(3)$ divides the conductor ideal $\mathfrak{f}(R)$.
Thus the upper bound in Corollary \ref{cor:chlthm631} is $18$.

Using Magma, we see that $\#\clb(R)=9$. In particular, $\#(\mathrm{Cl}(R)\backslash\overline{\mathrm{Cl}}(R))=3$. Thus neither of the bounds above is sharp.
\end{enumerate}
\end{example}

\begin{remark}\label{rmk:comparetwoex}
The two examples in Example \ref{ex:bounds} show that the relative sharpness of the two upper bounds depends on the example: the bound in Corollary \ref{cor:chlthm631} is sharper in Case (1), whereas the bound in Corollary \ref{thm:bassupperbound_global1} is sharper in Case (2).
\end{remark}

\subsection{The case that $R\cong \mathbb{Z}[x]/(x^3-mx^2+(m-1)x-1)$: an order associated with similarity classes of Cappell-Shaneson matrices}\label{subsec:Cap}
In this subsection, we provide an upper bound for $\#\overline{\mathrm{Cl}}(R)$ for an order with $m\in \mathbb{Z}$
$$R=\mathbb{Z}[\Theta_m]=\mathbb{Z}[x]/(x^3-mx^2+(m-1)x-1).$$  
As mentioned at the beginning of Section \ref{sec:application}, this is a Bass order in the cubic number field $\mathbb{Q}[x]/(x^3 - m x^2 + (m-1) x - 1)$ (see \cite[page 44]{AR84}).
The ideal class monoid $\clb(\mathbb{Z}[\Theta_m])$ plays a significant role in a topological problem concerning the \textit{Cappell-Shaneson homotopy 4-spheres}.

\subsubsection{\textbf{Background on Cappell-Shaneson matrices}}\label{subsub;cs}
A Cappell--Shaneson homotopy 4-sphere, originally constructed in \cite{Cap761}, is a smooth homotopy 4-sphere that admits an open book decomposition 
with fiber a punctured 3-torus. These manifolds remain important potential counterexamples to the smooth 4-dimensional Poincar\'e conjecture, 
which asserts that every smooth 4-manifold homotopy equivalent to $S^4$ is diffeomorphic to $S^4$.

Cappell--Shaneson homotopy 4-spheres are indexed by pairs $(\varepsilon, A)$, where $\varepsilon \in \mathbb{Z}/2\mathbb{Z}$ is a framing and $A \in \mathrm{SL}_3(\mathbb{Z})$ 
is a Cappell--Shaneson matrix satisfying $A - I \in \mathrm{SL}_3(\mathbb{Z})$. For a fixed framing, if two Cappell--Shaneson matrices are similar in $\mathrm{SL}_3(\mathbb{Z})$, 
then the corresponding spheres are diffeomorphic.  

The study of similarity classes of Cappell--Shaneson matrices is closely connected to the ideal class monoid 
$\overline{\mathrm{Cl}}(\mathbb{Z}[\Theta_m])$, since there exists a one-to-one correspondence between $\overline{\mathrm{Cl}}(\mathbb{Z}[\Theta_m])$ 
and the set of the similarity classes of Cappell--Shaneson matrices with trace $m$ 
(see \cite[page~44]{AR84}). Thus, to determine whether every Cappell--Shaneson homotopy 4-sphere is diffeomorphic to $S^4$, 
it suffices to analyze the spheres corresponding to each ideal class in $\overline{\mathrm{Cl}}(\mathbb{Z}[\Theta_m])$.

Consequently, understanding the size of $\overline{\mathrm{Cl}}(\mathbb{Z}[\Theta_m])$ is of considerable significance: 
the smaller the monoid, the fewer cases must be considered. In particular, establishing an upper bound for 
$\#\overline{\mathrm{Cl}}(\mathbb{Z}[\Theta_m])$ provides an effective estimate of the complexity of the diffeomorphism 
classification problem for Cappell--Shaneson spheres.

\subsubsection{\textbf{A refined upper bound for $\#\overline{\mathrm{Cl}}(R)$}}

We set up notations as follows:
\begin{itemize}
    \item For an integer $m\in \mathbb{Z}$,
    put $\phi(x)=x^3-mx^2+(m-1)x-1$ so that $R=\mathbb{Z}[x]/(\phi(x))$ and  $E=\mathbb{Q}[x]/(\phi(x))$.
   Note that $\phi(x)$ is irreducible over $\mathbb{Z}$.  
We emphasize that $\phi(x)$, $R$, and $E$ depend on the integer $m$, which we omit from the notation. 

    \item
We denote by $\Delta_{E_p}$ the discriminant ideal  $\Delta_{E_p/\mathbb{Q}_p}$.
Recall that $\ord_p(\Delta_{E_p})$ denotes the exponential valuation of a generator of $\Delta_{E_p}$ as an ideal of $\mathbb{Z}_p$.

\item In contrast to the previous notion, the discriminant of the number field $E$ over $\mathbb{Q}$ is a well-defined integer. We will denote it by  $\Delta_E  \left(\in \mathbb{Z}\backslash \{0\}\right)$ and by $|\Delta_E|$ its absolute value.

    \item The discriminant $\Delta_\phi$ of the polynomial $\phi(x)$ is $m^4-10m^3+31m^2-30m-23$. 
    Then
\begin{equation}\label{deltaphiposi}
    \textit{$m\geq 6$ or $m\leq -1$ $\Longleftrightarrow $ $\Delta_{\phi}>0\Longleftrightarrow$  $\Delta_E>0 \Longleftrightarrow  E$ is totally real.}
\end{equation}

    \item Define the constant $C_{\phi}=27\phi(m/3)=-2m^3+9m^2-9m-27$.

\end{itemize}

Our strategy to obtain an upper bound for $\#\overline{\mathrm{Cl}}(R)$ is to use Corollary \ref{thm:n=23} and Theorem \ref{thm:clbound}.
In order to use Corollary \ref{thm:n=23}, we will first investigate the Serre invariant at a prime integer $p$. 

Remark \ref{eq:about_serre}.(3) yields that the Serre invariant $S_p(R_p)=0$, equivalently $R_p=\mathcal{O}_{E_p}$, for all but finitely many prime $p$'s.
In particular, if $p\nmid \Delta_\phi$, then the reduction of $\phi(x)$ modulo $p$ has distinct roots, and hence $S_p(R_p)=0$. Therefore, it suffices to consider the case $p \mid \Delta_\phi$, in which the reduction of $\phi(x)$ modulo $p$ has a multiple root.

From now on until the end of this section, we suppose that $p \mid \Delta_\phi$. We distinguish the cases for a prime integer $p$ according to the type of a factorization of $\phi(x)$ in $\mathbb{Z}_p[x]$ as follows:
\begin{enumerate}
    \item[Case 1.] $\phi(x)$ is an irreducible polynomial in $\mathbb{Z}_p[x]$;

\item[Case 2.] $\phi(x)=\phi_{p,1}(x)\phi_{p,2}(x)$ where $\phi_{p,1}$ and $\phi_{p,2}$ are irreducible monic polynomials of degrees $2$ and $1$, respectively, and $R_p$ is a local ring;
\item[Case 3.] $\phi(x)=\phi_{p,1}(x)\phi_{p,2}(x)$ where $\phi_{p,1}$ and $\phi_{p,2}$ are irreducible monic polynomials of degrees $2$ and $1$, respectively, and $R_p$ is not a local ring;

\item[Case 4.] $\phi(x)$ is a product of three different monic polynomials of degree $1$ in $\mathbb{Z}_p[x]$.
\end{enumerate}
\begin{remark}\label{remark:pequals23}
    Notice that $2,3 \nmid \Delta_\phi$ for any $m$, so that $2$ and $3$ do not appear in Cases 1-4.

    Indeed, for $p=2$, evaluating $\overline{\Delta}_\phi=m^4+m^2+1$ at $m=0,1$ shows that $\Delta_\phi\not\equiv 0 \pmod{2}$, and hence $2\nmid \Delta_\phi$.
    Similarly, for $p=3$, evaluating $\overline{\Delta}_\phi=m^4+2m^3+m^2+1$ at $m=0,1,2$ shows that $\Delta_\phi\not\equiv 0 \pmod{3}$, and thus $3\nmid \Delta_\phi$. 

\end{remark}

\begin{lemma}\label{thm:Serreinv}
Let $p$ be a prime with $p\mid \Delta_\phi$.
If $R_p=\mathcal{O}_{E_p}$, then $S_p(R_p)=0$. 
If $R_p\neq \mathcal{O}_{E_p}$, then we have the following description for $E_p/\mathbb{Q}_p$ and $S_p(R_p)$:

\begin{itemize}
    \item In  Case 1, $E_p/\mathbb{Q}_p$ is ramified and $S_p(R_p)=1$.
    \item In Case 2, $E_p\cong E_{p,1}\oplus \mathbb{Q}_p$ where $E_{p,1}/\mathbb{Q}_p$ is ramified and $S_p(R_p)=1$.
    \item In  Case 3, $E_p\cong E_{p,1}\oplus \mathbb{Q}_p$ where $E_{p,1}/\mathbb{Q}_p$ is a quadratic field extension satisfying: 
    \begin{itemize}
        \item $E_{p,1}/\mathbb{Q}_p$ is unramified if and only if $2\mid \ord_p(\Delta_\phi)$. In this case $ S_p(R_p)=\frac{\ord_p(\Delta_\phi)}{2}$;
        \item $E_{p,1}/\mathbb{Q}_p$ is ramified if and only if $2\nmid \ord_p(\Delta_\phi)$. In this case $ S_p(R_p)=\frac{\ord_p(\Delta_\phi)-1}{2}$.
    \end{itemize}
    \item In  Case 4, $E_p\cong \mathbb{Q}_p^{\oplus3}$ and $S_p(R_p)=\frac{\ord_p(\Delta_\phi)}{2}$.
\end{itemize}
\end{lemma}
\begin{proof}
If $R_p=\mathcal{O}_{E_p}$, then $S_p(R_p)=0$ by the definition of the relative Serre invariant in Definition \ref{def:generalserreinv}. 
When $R_p\neq \mathcal{O}_{E_p}$, we treat Cases 1–4 separately.
\begin{enumerate}
    \item[Case 1.]
    Since $E_p$ is a cubic field extension of $\mathbb{Q}_p$ and  $R_p$ is a Bass order,  \cite[Proposition 3.5]{CHL} yields that  $E_p/\mathbb{Q}_p$ is totally ramified. 
    Then $S_p(R_p)=1$ by \cite[Proposition 3.9.(1)]{CHL}.

\item[Case 2.] 
By \cite[Proposition 6.8]{CHL}, after a suitable translation of the variable $x$ over $\mathbb{Z}_p$, the irreducible factors $\phi_{p,1}$ and $\phi_{p,2}$ of $\phi(x)$ are Eisenstein polynomials.
Then \cite[Lemma 6.10]{CHL} shows that $S_p(R_p)=\ord_{E_{p,1}}(\phi_{p,2}(\pi_{E_{p,1}})) =1$. 
Here, $\pi_{E_{p,1}}$ is a uniformizer of $E_{p,1}$ and 
$\ord_{E_{p,1}}(-)$ denotes the exponential valuation with respect to $\pi_{E_{p,1}}$.

\item[Case 3.] We write $E_p\cong E_{p,1}\oplus \mathbb{Q}_p$ where $E_{p,1}/\mathbb{Q}_p$ is a quadratic field extension.
By Remark \ref{remark:pequals23}, $p\neq 2,3$, so that  we have
\[
    \ord_p(\Delta_{E_p})=
    \left\{
    \begin{array}{l l}
    0 & \textit{if $E_{p,1}/\mathbb{Q}_p$ is unramified};\\
    1 & \textit{if $E_{p,1}/\mathbb{Q}_p$ is ramified}.
    \end{array}
    \right.
\]
Since $S_{p}(R_p)=\frac{\ord_p(\Delta_\phi)-\ord_p(\Delta_{E_p})}{2}$ is an integer by Equation (\ref{eq:serreinv}), it follows that $\ord_p(\Delta_\phi)$ is even if $E_{p,1}/\mathbb{Q}_p$ is unramified, and odd if $E_{p,1}/\mathbb{Q}_p$ is ramified.
Since $E_{p,1}/\mathbb{Q}_p$ is a quadratic field extension, it is either unramified or ramified. This yields the desired conclusion.

\item[Case 4.] It is clear that $E_p\cong \mathbb{Q}_p^{\oplus3}$. 
By \cite[Chapter I.3, Proposition 5.(iii)]{Cas}, we have $\ord_p(\Delta_{E_p})=0$.
This yields the desired conclusion.
\end{enumerate}
\end{proof}

In the following lemma, we will show that Cases 1–4 with $p \mid \Delta_\phi$ are characterized by the invariants $\Delta_\phi$ and $C_\phi$.

\begin{lemma}\label{prop:places}
    Let $p$ be a prime number such that $p\mid \Delta_\phi$.
    \begin{enumerate}
           \item 
        $p$ appears in either Case 1 or Case 2 if and only if $p\mid C_\phi$.
        In this case, we have $\mathrm{ord}_p(\Delta_\phi)= 2$, $3$, or $4$.
        More precisely, 
        \[
        \textit{$p$ appears in }\left\{
        \begin{array}{l l}
            Case \ 1 & \textit{if }\ord_p(\Delta_\phi)=2\textit{ or }4; \\
            Case \ 2 & \textit{if }\ord_p(\Delta_\phi)=3.
        \end{array}
        \right.
        \]
        In particular, $R_p=\mathcal{O}_{E_p}$ if and only if $\ord_p(\Delta_\phi)=2$.
        \item 
        Otherwise, that is $p\nmid C_\phi$ (equivalently $p$ appears in Case 3 or Case 4),
        \[
        \textit{$p$ appears in  }\left\{
        \begin{array}{l l}
            Case \ 3 & \textit{if }2\nmid \ord_p(\Delta_\phi) 
            \textit{ or }\left(\frac{\Delta_{\phi,\dot{p}}}{p}\right)=-1 ; \\
            Case \ 4 & \textit{if }2\mid \ord_p(\Delta_\phi) \textit{ and } \left(\frac{\Delta_{\phi,\dot{p}}}{p}\right)=1,
        \end{array}
        \right.
        \]where $\Delta_{\phi,\dot{p}}=\frac{\Delta_{\phi}}{p^{\ord_p(\Delta_\phi)}}$
         and $\left(\frac{\Delta_{\phi,\dot{p}}}{p}\right)$ denotes the Legendre symbol of $\Delta_{\phi,\dot{p}}$ with respect to $p$.
    \end{enumerate}
\end{lemma}
\begin{proof}
Recall that   $p\neq 2,3$ by Remark \ref{remark:pequals23}. 
Since $p\mid \Delta_\phi$, the polynomial $\overline{\phi}(x)$ has a multiple root in $\mathbb{F}_p$.
A prime number $p$ appears in Case 1 or Case 2 if and only if $R_p$ is a local ring, which is equivalent to $\overline{\phi}(x)$ having a triple root, by \cite[Section I.6, Lemma 4]{Ser79}.
Since $\phi''(x)=6x-2m$,  $\overline{\phi}(x)$ has a triple root if and only if $\overline{\phi}(m/3)=0$.
This is equivalent to the condition that $p \mid 27\phi(m/3)=C_{\phi}$.

\begin{enumerate}
    \item 
Suppose that $p\mid C_\phi$, equivalently $p$ appears in Case 1 or Case 2.
We first claim that $R_p=\mathcal{O}_{E_p}$ if and only if $\ord_p(\Delta_\phi)=2$, and in this case we claim that $p$ appears in Case 1.

Suppose that  $R_p=\mathcal{O}_{E_p}$. Since  $\mathcal{O}_{E_p}$ is a direct product of DVR's and $R_p$ is a local ring, $R_p=\mathcal{O}_{E_p}$ should be a DVR so that $E_p$ is a field extension of $\mathbb{Q}_p$. 
Between the two possibilities (Case 1 and Case 2), Case 2 is ruled out since $E_p$ is not a field in that case.
Hence, $p$ appears in Case 1. 

We now show that $\ord_p(\Delta_\phi)=2$.
The condition $p\mid \Delta_\phi$ implies that $E_p$ is totally ramified over $\mathbb{Q}_p$.
Thus we may and do assume that $\phi(x)$ is an Eisenstein polynomial by a suitable translation of a variable $x$ over $\mathbb{Z}_p$.
Indeed, the discriminant is invariant under this procedure. 
Since the discriminant of a monic cubic polynomial $x^3+bx^2+cx+d$ is 
$b^2c^2-4c^3-4b^3d-27d^2+18bcd$, we have that $\ord_p(\Delta_\phi)=\ord_p(-27d^2)=2$.

Conversely, suppose that $R_p\neq \mathcal{O}_{E_p}$. 
Then by a suitable translation of a variable $x$ over $\mathbb{Z}_p$, we may and do assume that $\overline{\phi}(x)=x^3$ in $\mathbb{F}_p[x]$.
In Case 1, $\ord_p(\phi(0))=2$  by \cite[Proposition 3.6]{CHL} since $R_p$ is a Bass order.
In Case 2, as in the proof of Lemma \ref{thm:Serreinv}, the irreducible factors $\phi_{p,1}(x)$ and $\phi_{p,2}(x)$ of $\phi(x)$ are Eisenstein polynomials  after a suitable translation of variable $x$ over $\mathbb{Z}_p$.
Thus, $\ord_p(\phi(0))=2$ so that $\ord_p(\Delta_\phi)>2$.

To complete the proof of (1), suppose that $\ord_p(\Delta_\phi)>2$, equivalently $R_p\neq \mathcal{O}_{E_p}$.
Since only Case 1 or Case 2 can occur, it suffices to show that $\ord_p(\Delta_\phi)=4$ if $p$ appears in Case 1, and $\ord_p(\Delta_\phi)=3$ if $p$ appears in Case 2.
By Lemma \ref{thm:Serreinv}, we have $S_p(R_p)=1$ when $p\mid C_\phi$ and $R_p\neq \mathcal{O}_{E_p}$. 
By Equation \eqref{eq:serreinv}, we have
$$\ord_p(\Delta_\phi)=\ord_p(\Delta_{E_p})+2S_p(R_p)=\ord_p(\Delta_{E_p})+2.$$
Thus it suffices to show that $\ord_p(\Delta_{E_p})=2$ if $p$ appears in Case 1, and $\ord_p(\Delta_{E_p})=1$ if $p$ appears in Case 2.

\begin{itemize}
    \item 
In Case 1, $E_p/\mathbb{Q}_p$ is ramified by Lemma \ref{thm:Serreinv}, and thus $\ord_p(\Delta_{E_p})=2$.

\item In Case 2, $E_p\cong E_{p,1}\oplus \mathbb{Q}_p$, where $E_{p,1}/\mathbb{Q}_p$ is ramified by Lemma \ref{thm:Serreinv}.
By \cite[Chapter I.3, Proposition 5.(iii)]{Cas}, $\ord_p(\Delta_{E_p})=\ord_p(\Delta_{E_{p,1}})=1$.
\end{itemize}

\item 
Suppose that $p\nmid C_\phi$, equivalently $p$ appears in Case 3 or Case 4.
In this case, we may write 
$\phi(x)=(x-a)(x^2+bx+c)$ over $\mathbb{Z}_p$.
We claim that $x^2+bx+c$ is irreducible if and only if $\Delta_\phi$ is not a square.
Suppose that the claim is true. 
If $\ord_p(\Delta_\phi)$ is odd or $\left(\frac{\Delta_{\phi,\dot{p}}}{p}\right)=-1$, then $\Delta_\phi$ is not a square and thus $p$ appears in Case 3. 
Otherwise, $\Delta_\phi$ is a square by Hensel's lemma, and thus $p$ appears in Case 4.

We now prove the claim. 
The polynomial $x^2+bx+c$ is irreducible if and only if its discriminant $b^2-4c$ is not a square in $\mathbb{Z}_p$.
    The claim  then follows from the formula:
$\Delta_\phi=(\lambda_1-\lambda_2)^2(\lambda_1-a)^2(\lambda_2-a)^2=(b^2-4c)(a^2+ab+c)^2$ 
for two roots $\lambda_1$ and $\lambda_2$ of $x^2+bx+c$.
\end{enumerate}
\end{proof}

\begin{corollary}\label{thm:boundmain}
We have the following upper bound for $\#\overline{\mathrm{Cl}}(R)$:
        \begin{equation}\label{eq:upperbound_of_main}
        \#\overline{\cl}(R)\leq \mathcal{A}(\Delta_\phi,C_\phi) \left(\frac{8}{3^7}\left(\frac{4}{\pi}\right)^{3r_2}|\Delta_E|^{\frac{3}{2}}+\frac{2}{3^4}\left(\frac{4}{\pi}\right)^{2r_2}|\Delta_E|+\frac{1}{3^3}\left(\frac{4}{\pi}\right)^{r_2}|\Delta_E|^{\frac{1}{2}}\right),
        \end{equation} 
        where \begin{align*}
        \mathcal{A}(\Delta_\phi&,C_\phi) =\underbrace{\left(\prod\limits_{\substack{ p\mid  C_{\phi}\\ \ord_p(\Delta_\phi)=4} }(p+1)\right)}_{\text{Case 1}}\underbrace{\left(\prod\limits_{\substack{p\mid  C_{\phi}\\ \ord_p(\Delta_\phi)=3 }}p\right)}_{\text{Case 2}}
      \underbrace{\left(\prod\limits_{\substack{p|\Delta_\phi,\ p\nmid  C_{\phi} \\2\nmid\ord_p(\Delta_\phi)}}
          \frac{p^{\frac{\ord_p(\Delta_\phi)+1}{2}}-1}{p-1}\right)}_{\text{Case 3 with $2\nmid \ord_p(\Delta_\phi)$}}\\ 
      &\underbrace{\left(\prod\limits_{\substack{p|\Delta_\phi,\ p\nmid C_{\phi} \\ 2\mid \ord_p(\Delta_\phi) \textit{ and }\left(\frac{\Delta_{\phi,\dot{p}}}{p}\right)=-1}}\frac{p^{\frac{\ord_p(\Delta_\phi)}{2}+1}+p^{\frac{\ord_p(\Delta_\phi)}{2}}-2}{p-1}\right)}_{\text{Case 3 with $2\mid \ord_p(\Delta_\phi)$}}
      \underbrace{\left(\prod\limits_{\substack{p|\Delta_\phi,\ p\nmid C_{\phi} \\ 2\mid \ord_p(\Delta_\phi)\textit{ and } \left(\frac{\Delta_{\phi,\dot{p}}}{p}\right)=1}}p^{\frac{\ord_p(\Delta_\phi)}{2}}\right)}_{\text{Case 4}}.
      \end{align*}
    Here, $\Delta_\phi=m^4-10m^3+31m^2-30m-23$, $C_{\phi}=-2m^3+9m^2-9m-27$, $\Delta_{\phi,\dot{p}}=\frac{\Delta_\phi}{p^{\ord_p(\Delta_\phi)}}$, and 
    $r_2=\left\{\begin{array}{l l}
         0& \textit{if } \Delta_E>0 \left(\textit{equivalently if } \Delta_\phi>0\right);\\
         1&  \textit{otherwise}.
    \end{array}\right.$
    \end{corollary}
\begin{proof}
We recall the notations from Section \ref{sec:application23}: $B_p(\phi)$ denotes an index set in bijection with irreducible factors $\phi_{p,i}$ of $\phi$ over $\mathbb{Z}_p$, and $R_{p,i}\cong \mathbb{Z}_p[x]/(\phi_{p,i}(x))$ for $i\in B_p(\phi)$.
Here, the invariant $\delta_p=\max_{i\in B_p(\phi)}\{S_p(R_{p,i})\}$, appearing in Corollary \ref{thm:n=23}, is computed as follows:
\begin{itemize}
    \item[Case 1.] 
    Since $\#B_p(\phi)=1$, we have $\delta_p=S_p(R_p)=1$ by Lemma \ref{thm:Serreinv}.
    \item[Case 2.]
    By the proof of Lemma \ref{thm:Serreinv}, the irreducible factors $\phi_{p,1}(x)$ and $\phi_{p,2}(x)$ of $\phi(x)$ are Eisenstein polynomials after a suitable translation of the variable $x$ over $\mathbb{Z}_p$.
    Therefore, both $R_{p,1}$ and $R_{p,2}$ are DVRs. 
    This yields that $\delta_p=0$.
    \item[Case 3.] 
    Since $R_p$ is not a local ring, the reductions $\overline{\phi_{p,1}}(x)$ and $\overline{\phi_{p,2}}(x)$ are coprime in $\mathbb{F}_p[x]$.
    Thus, it follows that $\phi_{p,1}(x)$ and $\phi_{p,2}(x)$ are coprime in $\mathbb{Z}_p[x]$. 
    By the Chinese remainder theorem, we have \[R_p \;\cong\; \mathbb{Z}_p[x]/(\phi_{p,1}(x)) \oplus \mathbb{Z}_p.\]
    This yields that $S_p(R_p)=S_p(\mathbb{Z}_p[x]/(\phi_{p,1}(x))) =\delta_p$.
    \item[Case 4.]  Since $R_{p,i}\cong\mathbb{Z}_p$ for each $i\in B_p(\phi)$, we have $\delta_p=0$.
\end{itemize}
By combining these computations with Lemma \ref{thm:Serreinv}, we obtain from Corollary \ref{thm:n=23} that
\[\#\overline{\mathrm{Cl}}(R)\leq \#\mathrm{Cl}(\mathcal{O}_E)\cdot \prod_{S_p(R_p)>0}\#(\Lambda_{E_p}\backslash X_{R_p}),
\]
where
  $$\#(\Lambda_{E_p}\backslash X_{R_p})=
    \left\{
\begin{array}{l l}
p+1 & \textit{in Case 1};\\
p & \textit{in Case 2};\\
1+(p+1)\frac{p^{\frac{\ord_p(\Delta_\phi)}{2}}-1}{p-1}& \textit{in Case 3 with $2\mid \ord_p(\Delta_\phi)$};\\
1+p\frac{p^{\frac{\ord(\Delta_\phi)-1}{2}}-1}{p-1}  & \textit{in Case 3 with $2\nmid \ord_p(\Delta_\phi)$};\\
p^{\frac{\ord_p(\Delta_\phi)}{2}} & \textit{in Case 4}.
\end{array}\right.$$    
Theorem \ref{thm:clbound} yields the following upper bound for $\#\mathrm{Cl}(\mathcal{O}_E)$
\[
\#\cl(\Mfo_E)\leq
\frac{\lfloor M \rfloor(\lfloor M \rfloor+1)(2\lfloor M \rfloor+1)}{6}\leq \frac{ M ( M +1)(2 M +1)}{6},
\]
where $M=\frac{2}{9}\left(\frac{4}{\pi}\right)^{r_2}|\Delta_E|^{1/2}$ in \eqref{eq:minkowski}.
Here, we recall that $\Delta_E>0$ if and only if $E$ is totally real (cf. \eqref{deltaphiposi}), and hence $r_2=\left\{\begin{array}{l l}
         0& \textit{if } \Delta_E>0;\\
         1&  \textit{otherwise}.
    \end{array}\right.$
This yields the desired conclusion.
\end{proof}

In (\ref{eq:upperbound_of_main}), we can rewrite $|\Delta_E|$ as a product formula using $\Delta_\phi$ and $C_\phi$.
\begin{proposition}\label{prop:delta_E}
We have
$$|\Delta_E|=\left(\prod\limits_{\substack{p\mid C_{\phi} \\ \ord_p(\Delta_\phi)=2 \textit{ or } 4}}p^2\right) \left(\prod\limits_{\substack{p\mid C_{\phi}\\ \ord_p(\Delta_\phi)=3 }} p\right)\left(\prod\limits_{\substack{p\mid \Delta_\phi,\ p\nmid C_{\phi} \\2 \nmid \ord_p(\Delta_\phi)}}p\right).$$

\end{proposition}

\begin{proof}
By the prime factorization of $|\Delta_E|\in\mathbb{Z}$, we have
\begin{equation}\label{eq:prodforDelta_E}
|\Delta_E|=\prod\limits_{p\in |\mathbb{Z}|} p^{\ord_p(\Delta_E)}.
\end{equation}
Since $\Delta_\phi=[\mathcal{O}_E:R]^2\Delta_E$ by \cite[Exercise 1.4.13]{Rei03}, we have $|\Delta_E|\mid |\Delta_\phi|$. Hence $\ord_p(\Delta_E)=0$ for $p\nmid \Delta_\phi$.
On the other hand, for $p\mid \Delta_\phi$, \cite[Chapter I.3, Proposition 6.(iii)]{Cas} shows that $\ord_p(\Delta_E)=\ord_p(\Delta_{E_p})$.
Putting these together, the Equation (\ref{eq:prodforDelta_E}) becomes
\[|\Delta_E|=\prod\limits_{p\mid \Delta_\phi}p^{\ord_p(\Delta_{E_p})}.\]

For each prime $p$ dividing $\Delta_\phi$, we compute $\ord_p(\Delta_{E_p})$ separately according to Cases 1–4.
From the proof of Lemma \ref{prop:places}.(1), we have 
$$\ord_p(\Delta_{E_p})=\left\{
\begin{array}{l l}
     2&  \textit{in Case 1};\\
     1&  \textit{in Case 2}.
\end{array}
\right.$$
In Case 3 or Case 4, since $S_p(R_p)=\frac{\ord_p(\Delta_\phi)-\ord_p(\Delta_{E_p})}{2}$ by Equation (\ref{eq:serreinv}), Lemma \ref{thm:Serreinv} yields that
$$\ord_p(\Delta_{E_p})=     \ord_p(\Delta_\phi)-2\lfloor \frac{\ord_p(\Delta_\phi)}{2} \rfloor=\left\{
\begin{array}{l l}
    1 & \textit{if $2\nmid \ord_p(\Delta_\phi)$}; \\
    0 & \textit{otherwise}.
\end{array}
\right.
$$
By combining these computations, we obtain the desired conclusion according to Lemma \ref{prop:places}.
\end{proof}

\begin{remark}\label{rmk:delta}
        For comparison with $|\Delta_E|$ in Proposition \ref{prop:delta_E}, we describe the prime factorization of $|\Delta_\phi|$ as follows:  
        \begin{equation}\label{eq:upperbound_Delta_phi}
        |\Delta_\phi|=\underbrace{\left(\prod\limits_{\substack{p\mid C_{\phi}\\\ord_p(\Delta_\phi)=2 }}p^2\right) \left(\prod\limits_{\substack{ p\mid C_{\phi}\\ \ord_p(\Delta_\phi)=4 }}p^4\right)}_{\text{Case 1}}
        \underbrace{\left(\prod\limits_{\substack{ p\mid C_{\phi}\\\ord_p(\Delta_\phi)=3}}p^3\right)}_{\text{Case 2}}
        \underbrace{\left(\prod\limits_{\substack{p\mid \Delta_{\phi}\\p\nmid C_\phi}}p^{\ord_p(\Delta_\phi)}\right)}_{\text{Case 3, 4}}.
        \end{equation}
        It is straightforward to see that $|\Delta_\phi|\geq |\Delta_E|$.
        In particular, for $p$ such that  $p\mid \Delta_\phi$ and $p\nmid C_\phi$, the exponent of $p$ in $|\Delta_\phi|$ may differ significantly from its exponent in $|\Delta_E|$.
\end{remark}
\begin{corollary}\label{cor:main_easy}
If  $\Delta_E>3075$, then we have the following upper bound for $\#\overline{\mathrm{Cl}}(R)$
\[
 \#\overline{\cl}(R) \leq \frac{2}{3^5} \Delta_\phi^{\frac{1}{2}}\cdot \Delta_E^{\frac{3}{2}}.
\]
\end{corollary}
\begin{proof}
We prove the argument using the upper bound for $\#\overline{\mathrm{Cl}}(R)$ given in Corollary \ref{thm:boundmain}.
Since $\Delta_E>3075$ (so that $r_2=0$),  the upper bound for $\#\cl(\Mfo_E)$ described in Theorem \ref{thm:clbound} is at most the following:

\[
 \left(
 \frac{8}{3^7}\Delta_E^{\frac{3}{2}}+\frac{2}{3^4}\Delta_E+\frac{1}{3^3}\Delta_E^{\frac{1}{2}}\right)\leq \frac{1}{3^5}\Delta_E^{\frac{3}{2}}.
\]

On the other hand, we recall that $p\mid \Delta_\phi$ only if $p>3$, by Remark \ref{remark:pequals23}.
By comparing $\mathcal{A}(\Delta_\phi,C_\phi)$ in Corollary \ref{thm:boundmain} with the product formula for $|\Delta_\phi|$ in (\ref{eq:upperbound_Delta_phi}) term by term,
we conclude
\[
\mathcal{A}(\Delta_\phi,C_\phi)\leq 2\Delta_\phi^{\frac{1}{2}}.
\]

\end{proof}

\begin{remark}\label{rmk:final}
In the above corollary, $\Delta_{E}>0$ if $m\geq 6$ or $m\leq -1$ by \eqref{deltaphiposi}. 
Therefore all but finitely many $E$'s satisfy the condition that $\Delta_E>3075$.

\end{remark}

\end{document}